% amstex file
\input amstex
\documentstyle {amsppt}
\magnification=\magstep1
\hsize=16truecm
\vsize=22.5truecm
\baselineskip=16truept
\NoRunningHeads
\NoBlackBoxes

\topmatter
\title Specialization of the torsion subgroup of the Chow group \endtitle
\author Chad Schoen \endauthor
\abstract An example is given in which specialization is not injective. \endabstract
\address Department of Mathematics, Duke University,
Box 90320, Durham, NC 27708-0320 USA,
e-mail: schoen\@ math.duke.edu . \endaddress
\thanks Partial support by the NSF (DMS-0200012) gratefully acknowledged \endthanks
\endtopmatter

\document

\define\op{\operatorname}
\define\zl{\Bbb Z_l}
\define\ql{\Bbb Q_l}
\define\inv{^{-1}}

\define\lbs{_{\bar s}}
\define\lbeta{_{\bar {\eta}}}
\define\lbsz{_{\bar s_0}}

Consider the diagram,
$$\CD
V_{\bar s}@>>> V_s@>>> \Cal V @<<< V_{\eta }@<<<V\lbeta \\
@VVV @VVV  @VVfV @VVV @VVV\\
\bar s@>\upsilon >>s@>i_s>> S @<g<< \eta @<\epsilon << \bar {\eta},
\endCD \tag 1$$
in which all squares are Cartesian, $S=Spec(A)$ is a finite type, regular, integral
affine $\Bbb Z$-scheme (or a localization of such), $g$ is the inclusion of
the generic point, $\eta =Spec(K)$,
$i_s$ is the inclusion of a non-generic point, $s=Spec(\Bbb F)$,
and $f$ is a smooth, projective morphism with geometrically
connected $n$-dimensional fibers. Furthermore, the map $\upsilon $ (respectively $\epsilon $)
corresponds to a choice of algebraic closure, $\Bbb F\subset \bar {\Bbb F}$ (respectively
$K\subset \bar K$).

Let $l$ be a prime number. For any abelian group, $B$, write $B[l^{\infty }]\subset B$ for
the subgroup consisting of elements annihilated by some power of $l$.
For $r\geq 0$ there is a specialization homomorphism of Chow groups,
$$\sigma ^r \lbs :CH^r(V\lbeta)[l^{\infty }]@>>>CH^r(V\lbs )[l^{\infty }].\tag 2$$
When $l$ is distinct from the characteristic of $\Bbb F$, $\sigma ^1\lbs $ is injective.
This fact has proved to be very useful for bounding torsion in
$CH^1(V_{\eta} )\subset CH^1(V\lbeta )^{Gal(\bar K/K)}$ (cf.
\cite {Sil, VIII.7.3.2}). It is interesting to ask if injectivity continues to hold for
$r>1$. In fact it does for $r\in \{2,n\}$
as will be recalled in the Proposition below. The purpose of this note is to
show that injectivity may fail in the range $2<r<n$.\footnote{After this paper was written,
C. Soul\'e and C. Voisin posted {\it Torsion cohomology classes and algebraic cycles on
complex projective manifolds} AG/0403254 on the e-print archives lanl.arXiv.org. Their
e-print is also concerned with non-injectivity of specialization of torsion in the
Chow group among other issues. Both the differences and the similarities between the
cycles considered in the e-print and the cycles
considered here appear to be interesting.} The failure
of injectivity when $l=char(\Bbb F)$ and $r$ is arbitrary is classical \cite {Sil, III.6.4}.

Before proving the injectivity and non-injectivity assertions, we describe the map,
$\sigma ^r_{\bar s}$, in detail. Extra effort is required here to deal with
the case $codim_S(s)>1$, which plays an important role in later arguments.

\subhead Construction of the specialization homomorphism \endsubhead
By localizing we may assume that $S$ is the spectrum of a regular local ring and that
$s$ is the closed point. Write $\tilde S,E@>>>S,s$ for the blow up of $S$ along $s$.
The local ring at the exceptional divisor, $\Cal O_{\tilde S,E}$, is the
valuation ring of a discrete valuation, $v:K^*\to \Bbb Z$. Choose a valuation
$\bar v$ of $\bar K$ extending $v$ \cite {Bou,VI.3.3}.
In each intermediate field, $K\subset L\subset \bar K$,
finite over $K$, $\bar v$ specifies a discrete valuation ring, $\Cal O_{L,\bar v}$,
finite over  $\Cal O_{\tilde S,E}$.
Write $s_{L,\bar v}$ for the spectrum of its residue field. Base change
$f$ by $Spec(\Cal O_{L,\bar v}) \to S$. There is a specialization homomorphism
\cite {Fu, 20.3.1},
$$\sigma ^r _{L,\bar v}:CH^r(V_L)\to CH^r(V_{s_{L,\bar v}}).$$
Let $L\subset L'$ be a finite extension of intermediate fields. Since
specialization is functorial for flat morphisms \cite {Fu, Proposition 20.3(b)},
the maps $\sigma ^r _{L,\bar v}$ and $\sigma ^r _{L',\bar v}$ are related by a
commutative diagram. Taking the limit over finite intermediate fields
gives a map,
$$\sigma ^r _{\bar K,\bar v}:CH^r(V\lbeta )\to CH^r(V_{s_{\bar K,\bar v}}),\tag 3$$
where $s_{\bar K,\bar v}$ is the spectrum of an algebraic closure of the residue field
of $\Cal O_{\tilde S,E}$. Now (2) is constructed from (3) by restricting to $l$-power
torsion subgroups and composing with the inverse of the pullback isomorphism \cite {Le},
$$CH^r(V_{\bar s})[l^{\infty }]@>\sim >> CH^r(V_{s_{\bar K,\bar v}})[l^{\infty }].$$
The construction is independent of the choice of $\bar v$: Any valuation on $\bar K$
extending $v$ has the form $\bar v\circ \gamma $ for $\gamma \in Gal(\bar K/K)$ \cite {Bou, VI.8.6}.
There is a commutative diagram,
$$\CD
CH^r(V\lbeta )[l^{\infty }]@>\sigma ^r _{\bar K,\bar v}>> CH^r(V_{s_{\bar K,\bar v}})[l^{\infty }]@<\sim <<
CH^r(V\lbs )[l^{\infty }]\\
@VVV @VVV @VVV \\
CH^r(V\lbeta )[l^{\infty }]@>\sigma ^r _{\bar K,\bar v\circ \gamma }>>
CH^r(V_{s_{\bar K,\bar v\circ \gamma }})[l^{\infty }]@<\sim << CH^r(V\lbs )[l^{\infty }],
\endCD $$
with vertical maps isomorphisms induced by $\gamma $.

\proclaim {Lemma 1} Let $\Cal X=\sum n_i\Cal X_i\in Z^r(\Cal V)$ be a linear combination of integral
subschemes which are
flat over $S$. Write $X_i$ for the generic fiber of $\Cal X_i$. Suppose that the image of $\sum n_iX_i$
in $CH^r(V\lbeta )$ is annihilated by a power of $l$. Then $\sigma ^r _{\bar s}(\sum n_iX_i)$ is equal to
the class of $\sum n_i(\Cal X_i\times _S\bar s)$ in $CH^r(V_{\bar s})$.
\endproclaim
\demo{Proof} This follows from standard intersection
theory \cite {Fu, 1.5, 6.1, 20.3}
\hfill $\square $\enddemo

\subhead An injectivity result \endsubhead
To place the main result of this note in perspective we recall

\proclaim{Proposition}
In the situation of \rom{(1)} suppose that $l\in A^*$. If $r\in \{ 1,2,n\} $, then
$\sigma ^r _{\bar s}$ is injective.
\endproclaim
\demo{Proof}
For a smooth, projective variety, $T$, over an algebraically closed field of
characteristic prime to $l$, Bloch has defined a cycle class map to etale
cohomology \cite {Bl},
$$\lambda ^r_{T}:CH^r(T)[l^{\infty }]\to H^{2r-1}(T,\ql /\zl (r)).$$
When $r=1$ this map is an isomorphism \cite {Bl, 3.6}. Roitman's theorem
\cite {Ro}, \cite {Bl, 4.2} shows that $\lambda ^n_T$ is an isomorphism when $n=dim(T)$.
As a consequence of the Merkuriev-Suslin theorem, $\lambda ^2_T$ is injective
\cite {Co-Sa-So, Corollaire 4}.

Since $f$ is smooth and projective,
the cospecialization map on cohomology,
$$c_s:H^{2r-1}(V\lbs ,\ql /\zl (r))\to H^{2r-1}(V\lbeta ,\ql /\zl (r)),$$
is an isomorphism \cite {Mi,VI.4.2}. The assertion follows from the
compatiblity of Bloch's cycle class map with specialization:
$c_s\inv \circ \lambda ^r_{V\lbeta }=\lambda ^r_{V\lbs }\circ \sigma ^r \lbs $
\cite {Bl, 3.8}.
\hfill $\square $\enddemo

\subhead The main result \endsubhead
\proclaim {Theorem}
There exist diagrams \rom{(1)} with $S$ of finite type over $Spec(\Bbb Z)$
and primes $l\in A^*$ such that for all $r$ in the range,
$2<r<n$, and all closed points $s\in S$, $\sigma ^r _{\bar s}$ is
not injective.
\endproclaim
\demo{Proof} Clearly we need a variety $V\lbeta $ for which $Ker(\lambda ^r_{V\lbeta })\neq 0$
for $r$ in the range, $2<r<n$. In \cite {Sch} certain varieties with this property
were constructed as products of a smooth projective variety, $W\lbeta$, and an elliptic
curve, $Y_{\lbeta}$. Cycles in $Ker(\lambda ^r)$ were constructed as exterior products, $z\times \tau $,
with $z\in Z^{r-1}(W\lbeta )$ and $\tau \in CH^1(Y\lbeta )[l^{\infty }]$. The idea of the proof
is to arrange that $z$ specialize to a torsion class and then to exploit the divisibility of
$\tau $ to conclude that $z\times \tau $ specializes to zero.
To implement this idea we begin with a presentation of the exterior product construction
of \cite {Sch, \S 1} in a relative context.

Fix a prime number $l$ and a finitely generated field, $k_0$, with $l\neq char(k_0)$.
Let $k_0\subset K$ be a finitely generated extension of transcendence degree one.
Let $Y/K$ be an elliptic curve whose $j$-invariant, $J(Y)\in K$, is transcendental
over $k_0$. Suppose that there is a torsion point, $\zeta _K\in Y$, of exact
order $l^m$. Let $W$ be a geometrically connected, smooth, projective variety over
$k_0$ of dimension $n-1$. Let $A_0\subset k_0$ be a regular integral domain of finite type over
$\Bbb Z$ with fraction field $k_0$. By inverting an element in
$A_0$ if necessary, we may arrange that
$l\in A_0^*$ and that $W$ extends to smooth, projective morphism,
$h:\Cal W\to S_0:=Spec(A_0)$.
Let $A$ be a regular integral domain, flat and of finite type over $A_0$ whose fraction
field may be identified with $K$. By inverting an element in $A$ if necessary, we
arrange that $Y/K$ extends to an abelian scheme, $\Cal Y\to S:=Spec(A)$. The
identity section is denoted $e$ and the $l^m$-torsion section which extends
$\zeta _K$, $\zeta $.
The composition
$$f\ :\ \Cal V:=\ \Cal W\times _{S_0}\Cal Y\ @>>>\ \Cal Y\ @>>>\ S$$
is projective and smooth of relative dimension $n$ with connected fibers.
The fibers have product structures:
$$\align
V_{\eta}=&\ \Cal W\times _{S_0}\Cal Y\times _S\eta \simeq \Cal W\times _{S_0}\Cal Y_{\eta }
\simeq W\times _{k_0}\Cal Y_{\eta } \\
V_s=&\ \Cal W\times _{S_0}\Cal Y\times _Ss \simeq \Cal W\times _{S_0}\Cal Y_s
\simeq \Cal W_{s_0}\times _{s_0}\Cal Y_s,
\endalign $$
where $s_0\in S_0$ is the image of $s$ and $s$ is an arbitrary closed point of $S$.

The elements of $CH^r(V\lbeta )[l^{\infty }]$ whose specialization we wish to study are
constructed from subschemes of $\Cal Y$ (respectively $\Cal W$) which are flat over
$S$ (respectively $S_0$). Define groups of algebraic cycles
$$\align
Z^{r-1}_{fl}(\Cal W)&=\{ \sum n_i\Cal Z_i\in  Z^{r-1}(\Cal W): \text{each subscheme}\  \Cal Z_i
\ \text{is flat over}\  S_0\}.\\
Z^1_{fl}(\Cal Y)&=\{ \sum n_i\Cal U_i\in  Z^1(\Cal Y): \text{each subscheme}\  \Cal U_i
\ \text{is flat over}\  S\}.\\
Z^r_{fl}(\Cal V)&=\{ \sum n_i\Cal X_i\in  Z^r(\Cal V): \text{each subscheme}\  \Cal X_i
\ \text{is flat over}\  S\}.
\endalign $$
The image of $Z^1_{fl}(\Cal Y)\to CH^1(\Cal Y)$ is denoted $CH^1_{fl}(\Cal Y)$.
Define $CH^{r-1}_{fl}(\Cal W)$ and $CH^r_{fl}(\Cal V)$ analogously.
Let $\bar k_0$ be an algebraic closure of $k_0$. We construct a diagram
$$\CD
CH^{r-1}(\Cal W_{\bar s_0})\otimes CH^1(\Cal Y\lbs )@>\times >>
CH^r(\Cal W_{s_0}\times _{s_0}\Cal Y\lbs ) \\
@AAi^{!}_W\otimes i^!_YA @AAi^!_{V}A \\
CH^{r-1}_{fl}(\Cal W)\otimes CH^1_{fl}(\Cal Y)@>\times >>
CH^r_{fl}(\Cal W\times _{S_0}\Cal Y) \\
@VVj_W^*\otimes j_Y^*V @VVj^*_{V}V \\
CH^{r-1}(W_{\bar k_0})\otimes CH^1(\Cal Y\lbeta )@>\times >>
CH^r(\Cal W_{\bar k_0}\times _{\bar k_0}\Cal Y\lbeta ).
\endCD $$
The horizontal maps are exterior product maps \cite {Fu, 1.10 and proof of Proposition 20.2}
\cite {Sch, 1.1}. In the middle row an element represented by $\Cal Z\otimes \Cal U$ for
subschemes $\Cal Z\subset \Cal W$ flat over $S_0$ and $\Cal U\subset \Cal Y$ flat over $S$ is
mapped to the class of the subscheme $\Cal Z\times _{S_0} \Cal U$, which is flat over
$S$ \cite {Ha,III.9.2b,c}.
The maps labeled $j^*$ are flat pullback maps, while the maps $i^!$ are intersections
with a geometric closed fiber. All of the maps in the diagram
may be defined on the level of algebraic cycles (ie. without first passing
to rational equivalence classes). The top square of the diagram commutes. In
fact with $\Cal Z$ and $\Cal U$ as above,
$$i^!_V(\Cal Z\otimes \Cal U)=\Cal Z\times _{S_0}\Cal U\times _S\bar s \simeq
\Cal Z_{s_0}\times _{s_0}\Cal U_{\bar s} \simeq \Cal Z_{\bar s}\times _{\bar s}\Cal U_{\bar s}
=i^!_W(\Cal Z)\times i^!_Y(\Cal U). $$
The bottom square also commutes:
$$j^*_V(\Cal Z\times _{S_0}\Cal U)=\Cal Z\times _{S_0}\Cal U\times _S\bar {\eta }
\simeq \Cal Z_{k_0}\times _{k_0}\Cal U_{\bar {\eta }}
\simeq \Cal Z_{\bar k_0}\times _{\bar k_0}\Cal U_{\bar {\eta }}
=j^*_W(\Cal Z)\times j^*_Y(\Cal U).$$

Define $\Cal T=\zeta -e\in CH^1_{fl}(\Cal Y)$.

\proclaim{Lemma 2}
(i) If $\Cal Z\in CH^{r-1}_{fl}(\Cal W)$ is such that $j^*_W(\Cal Z)\otimes 1/l^m$
has exact order $l^m$
in $CH^{r-1}(W_{\bar k_0})\otimes \ql /\zl $, then $j^*_W(\Cal Z)\times j^*_Y(\Cal T)$
has exact order $l^m$ in $CH^r(V\lbeta )[l^{\infty }]$.

(ii) If $i^!_W(\Cal Z)\in CH^{r-1}(\Cal W_{\bar s_0})_{tors}$, then
$\sigma ^r \lbs (j^*_W(\Cal Z)\times j^*_Y(\Cal T))=0$ in $CH^r(V\lbs )[l^{\infty }]$.
\endproclaim
\demo{Proof}
(i) There is an injective group homomorphism, $\ql /\zl \to CH^1(\Cal Y_{\bar {\eta }})$,
mapping $1/l^m$ to $j^*_Y(\Cal T)$, whose image is a direct summand of
$CH^1(\Cal Y_{\bar {\eta }})$. Thus $j^*_W(\Cal Z)\otimes j^*_Y(\Cal T)\in
CH^{r-1}(W_{\bar k_0})\otimes CH^1(\Cal Y_{\bar {\eta }})$ has exact order
$l^m$. By \cite {Sch, 0.2} the exterior product map,
$$CH^{r-1}(W_{\bar k_0})\otimes CH^1(\Cal Y\lbeta )[l^{\infty}]@>\times >>
CH^r(\Cal W_{\bar k_0}\times _{\bar k_0}\Cal Y\lbeta )[l^{\infty }],$$
is injective.

(ii) $j^*_W(\Cal Z)\times j^*_Y(\Cal T)=j^*_V(\Cal Z\times \Cal T)$ and
$\sigma ^r \lbs (j^*_V(\Cal Z\times \Cal T))=i^!_V(\Cal Z\times \Cal T)$ by Lemma 1.
Observe that $i^!_W(\Cal Z)\otimes i^!_Y(\Cal T)=0\in CH^{r-1}(\Cal W_{\bar s})
\otimes CH^1(\Cal Y\lbs )[l^{\infty }]$, since $CH^1(\Cal Y\lbs )[l^{\infty }]$
is a divisible group and $i^!_W(\Cal Z)$ is torsion. Now the commutativity of
the upper square in the diagram implies $\sigma ^r \lbs (j^*_V(\Cal Z\times \Cal T))=0$.
\hfill $\square $\enddemo

To complete the proof of the theorem we give an example where the hypotheses of
Lemma 2 are fulfilled.

\medpagebreak
\noindent {\bf Jacobians of genus three curves.}
Set $k_0=\Bbb Q(a)$,
$$E:\ \ y^2z-[(a^2-4)x^3+(2a^2-4a)x^2z+(a^2-4)xz^2]=0,$$
and $W=E^3$. Now $W$ is isogenous to the jacobian of the genus $3$ curve,
$$C:\ \ x^4+y^4+z^4+a(x^2y^2+y^2z^2+z^2x^2)=0.$$
A map, $\pi :C\to E,$ is given by the field extension, $k_0(E)\to k_0(C)$ \cite {Bu-Sch-Top,4.2},
$$x/z \mapsto (2(x/z)^2+(y/z)^2+a)^2,\qquad
y/z\mapsto (y/z)(2(x/z)^2+(y/z)^2+a).$$
The element, $\kappa \in \op{Aut}(\Bbb P^2)$, given by
$$x\circ \kappa =y,\qquad y\circ \kappa =z,\qquad z\circ \kappa =x,$$
stabilizes $C$ and gives rise to an embedding,
$$\varrho :C@>>>E^3, \qquad \varrho =
(\pi ,\pi \circ \kappa ,\pi \circ \kappa ^2).$$
Define a cycle,
$$Z :=\varrho (C)-(-1)_*\varrho (C)\in Z^2(W),$$
where $(-1)$ denotes inversion in the group law on the abelian variety, $W$.
For certain values of $a$ and $l$ calculations based on \cite {Bl-Es},
and \cite {Bu-Sch-Top} show that the class of $Z$ in $CH^2(W_{\bar k_0})$ is
not torsion and its image in $CH^2(W_{\bar k_0})/l$ is not zero. For instance
this holds when $a=l=5$ \cite {Sch2, 3.1}, \cite {Bu-Sch-Top, 4.8, 9.2} and when
$a$ is an indeterminant and $l\in \{ 5,7,11,13,17\} $ \cite {Sch2, 3.3}.
Choose $A_0\subset k_0$ as above, so that $C$ and $E$ extend to relative curves,
smooth and projective over $S_0=Spec(A_0)$, and $\rho $ extends as well. Write
$\Cal Z\in Z^2_{fl}(\Cal W)$ for the obvious extension of $Z$.
Now $j^*_W(\Cal Z)\otimes 1/l^m=Z\otimes 1/l^m$ has exact order $l^m$
in $CH^2(W_{\bar k_0})\otimes \ql /\zl $.

Consider $i^!_W(\Cal Z)=\rho \lbsz (C\lbsz )-(-1)_*\rho \lbsz (C\lbsz )\in CH^2(W\lbsz )$, where
$C\lbsz $ is the fiber of the relative curve extending $C$ over $\bar s_0$. The cohomology class,
$[i^!_W(\Cal Z)]\in H^4(W\lbsz ,\ql (2))$, vanishes, since inversion on the abelian variety, $W\lbsz $,
acts trivially on $H^4(W\lbsz ,\ql (2))$.  Since $s_0\in S_0$ is a closed point,
$\bar s_0$ is the spectrum of an algebraic closure of a finite field. Soul\'e's theorem
\cite {So, Th\'eor\`eme 4}
says that nullhomologous, one dimensional cycles on an Abelian variety over an
algebraic closure of a finite field are torsion, ie. $i^!_W(\Cal Z)\in CH^2(W_{\bar s_0})_{tors}$.
Thus the hypotheses of Lemma 2 are satisfied.

This proves the theorem when $n=4$. It is easy to extend this to the case
$n>4$ by taking $W=E^3\times \Bbb P^{n-4}$ instead of $W=E^3$ (cf. \cite {Sch, 3.3}).
\hfill $\square $\enddemo

\example{Remark}
The arguments in this paper do not yield results about the specialization map (2) when
$S$ is an open subset of the ring of integers of a number field. Indeed the assumption
that $J(Y)\in K$ is {\it transcendental} over $k_0$ is needed here because of the important
role that this hypothesis plays in \cite {Sch, 0.2}.
Totaro \cite {To, 7.2} has contructed elements in the kernel of Bloch's cycle class map
which are defined over number fields.
It would be interesting to know how these cycles behave under specialization.
\endexample

\subhead References \endsubhead

[Bl] Bloch, S., {\it Torsion algebraic cycles and a theorem of Roitman},
Compositio Math. {\bf 39}, (1979), 107-127.

[Bl-Es] Bloch, S. and Esnault, H., {\it The coniveau filtration and non-divisibility
for algebraic cycles}, Math. Ann. 304, p. 303-314 (1996).

[Bou] Bourbaki, N., Commutative Algebra (Chapters 1-7), Springer-Verlag (1989)

[Bu-Sch-Top] Buhler, J., Schoen, C., and Top, J., {\it Cycles, $L$-functions
and productes of elliptic curves}, J. reine angew. Math. 492, p. 93-133 (1997)

[Co-Sa-So] Colliot-Th\'el\`ene, J.-L., Sansuc, J.-J., and Soul\'e, C.,
{\it Torsion dans le groupe de Chow de codimension deux}, Duke Math. J. {\bf 50} (1983),
763--801.

[Fu] Fulton, W., Intersection Theory. Springer-Verlag, New-York, (1984).

[Ha] Hartshorne, R., Algebraic Geometry, Springer-Verlag, New York, (1977).

[Le] Lecomte, F., {\it Rigidit\'e des groupes de Chow}, Duke Math. J. 53, p. 1011-1046 (1986).

[Mi] Milne, J.S., \'Etale Cohomology, Princeton University Press, Princeton, (1980).

[Ro] Rojtman, A. A., {\it The torsion of the group of $0$-cycles modulo rational equivalence},
Ann. of Math. (2)  111, p. 553-569  (1980)

[Sch] Schoen, C., {\it On certain exterior product maps of Chow groups,
 Mathematical Research Letters}, {\bf 7}, p. 177-194 (2000).

[Sch2] Schoen, C., {\it The Chow group modulo l for the triple product of a general elliptic
curve}, Asian Journal of Mathematics {\bf 4}, 987-996, (2000).

[Sil] Silverman, J., The Arithmetic of Elliptic Curves, Springer-Verlag, New York, (1986).

[So] Soul\'e ,C., {\it Groupes de Chow et $K$-th\'eorie de vari\'et\'es sur un corps
fini}, Math.\ Ann.\ {\bf 268} (1984), 317-345.

[To] Totaro, B., {\it Torsion algebraic cycles and complex cobordism},
J. Amer. Math. Soc.  {\bf 10}  (1997),  467-493.
\enddocument